\documentclass[11pt,a4paper]{article}

\usepackage{url,amsmath,amssymb,latexsym,pstricks,mathrsfs,comment,amsthm,graphicx,tikz,tikz-cd,enumerate,accents,pgffor,cite,wrapfig,multicol,float,cases,calc,bibspacing,geometry,bbm,arydshln}
\usepackage[colorlinks]{hyperref}
\usepackage[T1]{fontenc}
\usepackage{ wasysym }
\usepackage{stmaryrd}
\usepackage{enumitem}

\usepackage{hhline}

\newcounter{ncols}
\newcounter{incols}

\usetikzlibrary{calc}

\allowdisplaybreaks

\numberwithin{equation}{section}

\newtheorem{thm}[equation]{Theorem}
\newtheorem{lemma}[equation]{Lemma}

\theoremstyle{definition}

\newcommand{\nc}{\newcommand}
\nc{\rnc}{\renewcommand}

\let\oldproofname=\proofname
\rnc{\proofname}{\rm\bf{\oldproofname}}

\nc{\pf}{\begin{proof}}
\nc{\epf}{\end{proof}}

\rnc{\arraystretch}{1.2}

\nc\ord{\operatorname{ord}}

\nc\id{\operatorname{id}}
\nc\ba{{\bf a}}
\nc\bb{{\bf b}}
\nc\A{\mathcal A}
\rnc\S{\mathcal S}
\nc\Z{\mathbb Z}
\nc\bone{{\bf 1}}

\nc\ben{\begin{enumerate}[label=\textup{(\roman*)},leftmargin=7mm]}
\nc\een{\end{enumerate}}
\nc{\bit}{\begin{itemize}}
\nc{\eit}{\end{itemize}}
\nc{\firstpfitem}[1]{#1.}
\nc{\pfitem}[1]{\medskip \noindent #1.}

\nc\T{\mathcal T}

\nc{\al}{\alpha}
\nc{\be}{\beta}
\nc{\ga}{\gamma}
\nc{\de}{\delta}

\nc{\rank}{d}

\nc{\set}[2]{\{ {#1} : {#2} \}} 
\nc{\sub}{\subseteq}
\nc{\la}{\langle}
\nc{\ra}{\rangle}

\nc{\COMMA}{,\qquad}
\nc{\AND}{\qquad\text{and}\qquad}
\nc{\ANDSIM}{\qquad\text{and similarly}\qquad}

\nc{\lv}[1]{\fill (#1,0)circle(.09);}
\nc{\lvs}[1]{{\foreach \x in {#1} { \lv{\x}}}}
\nc{\stline}[2]{\draw(#1,2)--(#2,0);}
\nc{\stlines}[1]{{\foreach \x/\y in {#1} { \stline{\x}{\y} }}}

\begin{document}

\title{\vspace{-1cm} Generating wreath products of symmetric and alternating groups}
\author{James East\footnote{Centre for Research in Mathematics and Data Science, Western Sydney University, Locked Bag 1797, Penrith NSW 2751, Australia. {\it Email:} {\tt j.east\,@\,westernsydney.edu.au}.  Supported by ARC Future Fellowship FT190100632.}
\ and
James Mitchell\footnote{Mathematical Institute, School of Mathematics and Statistics, University of St Andrews, St Andrews, Fife KY16 9SS, UK. {\it Email:} {\tt jdm3\,@\,st-andrews.ac.uk}.}}
\date{}

\maketitle

\vspace{-1cm}

\begin{abstract}
We show that the wreath product of two finite symmetric or alternating groups is 2-generated.


\end{abstract}


\section{Introduction and statement of the main result}

For a positive integer $n$, we write $\S_n$ and $\A_n$ for the symmetric and alternating groups of degree $n$, regarded as groups of permutations of $\{1,\ldots,n\}$ in the usual way.  We write permutations to the right of their arguments, and compose from left to right.

\begin{figure}[h]
\begin{center}
\scalebox{.9}{
\begin{tikzpicture}[scale=1]
\tikzstyle{vertex}=[circle,draw=black, fill=white, inner sep = 0.06cm]
\draw[-{latex}](6+1,-1)--(10-1,-1);
\begin{scope}[shift={(0,0)}]	
\lvs{1,...,5}
\stlines{1/2,2/3,3/4,4/1,5/5}
\draw(.5,1)node[left]{$(\ba,f)=$};
\foreach \x in {1,...,5} {\node[vertex] () at (\x,2){\small $a_{\x}$};}
\foreach \x in {1,...,5} {\draw[dotted](\x,0)--(\x,-2);}
\end{scope}
\begin{scope}[shift={(0,-4)}]	
\lvs{1,...,5}
\stlines{1/2,2/1,3/4,4/5,5/3}
\draw(.5,1)node[left]{$(\bb,g)=$};
\foreach \x in {1,...,5} {\node[vertex] () at (\x,2){\small $b_{\x}$};}
\end{scope}
\begin{scope}[shift={(10,-2)}]	
\lvs{1,...,5}
\stlines{1/1,2/4,3/5,4/2,5/3}
\draw(5.5,1)node[right]{$=(\ba,f)(\bb,g)$};
\foreach \x/\y in {1/2,2/3,3/4,4/1,5/5} {\node[vertex] () at (\x,2){{\scriptsize $a_{\x}b_{\y}$}};}
\end{scope}
\end{tikzpicture}
}
\end{center}
\vspace{-5mm}
\caption{Calculating a product in $G\wr S$.}
\label{fig:prod}
\end{figure}
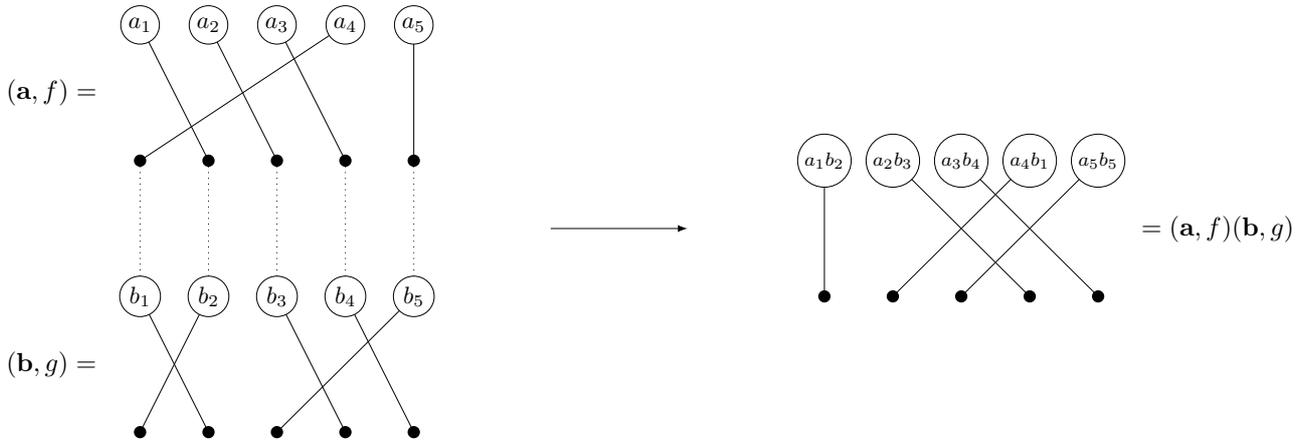

For a subgroup $S$ of $\S_n$, and an arbitrary group $G$, the wreath product~${G\wr S}$ is the semidirect product $G^n\rtimes S$, where $S$ acts on $G^n$ by permuting coordinates.  Formally, if ${\ba=(a_1,\ldots,a_n)\in G^n}$, and if $f\in S$, then we define $f\ba=(a_{1f},\ldots,a_{nf})$.  The elements of $G\wr S$ are of the form $(\ba,f)$, where $\ba\in G^n$ and $f\in S$, and such a pair will often be denoted by $(\ba,f)=(a_1,\ldots,a_n;f)$.  The product in $G\wr S$ is given by $(\ba,f)(\bb,g)=(\ba(f\bb),fg)$, or in the alternative form by
%
\[
(a_1,\ldots,a_n;f)(b_1,\ldots,b_n;g) = (a_1b_{1f},\ldots,a_nb_{nf};fg).
\]
There is a useful diagrammatic interpretation of this product, which will help with many of the calculations below.  Figure \ref{fig:prod} gives an example in the case that $n=5$, and $f=(1,2,3,4)$ and $g=(1,2)(3,4,5)$.

As usual, for a group $G$, we write $\rank(G)$ for the minimum number of elements required to generate~$G$.  Throughout the paper, generation always means generation as a semigroup: i.e., every element of the group is a product of positive powers of the generators.  For a finite non-trivial group, this is of course equivalent to generation as a group (where negative powers of the generators are allowed), but some later results apply to finitely-generated infinite groups.
%
Here is our main result: 

\begin{thm}\label{thm:main}
If $G$ and $S$ are arbitrary finite symmetric or alternating groups, then $\rank(G\wr S)\leq2$.  Moreover, $\rank(G\wr S)=1$ if and only if one of the following holds:
\ben
\item \label{cyc1} $G$ is one of $\S_1({}=\A_1)$ or $\A_2$, and $S$ is one of $\S_1$, $\S_2$, $\A_2$ or $\A_3$, or
\item \label{cyc2} $G$ is one of $\S_1$, $\S_2$, $\A_2$ or $\A_3$, and $S=\S_1$.
\een
\end{thm}

Our interest in Theorem \ref{thm:main} comes primarily from the paper \cite{AS2009}, the main result of which determines the minimum number of generators for a semigroup of transformations preserving a uniform partition of a finite set.  Such a semigroup can be thought of as a wreath product $\T_m\wr\T_n$ of two full transformation semigroups, and a key step was to show that a wreath product $\S_m\wr\S_n$ is 2-generated, which is of course a special case of Theorem \ref{thm:main}.  This fact about $\S_m\wr\S_n$ was proved in \cite{AS2009} using rather sophisticated methods from representation theory; see also \cite{BLS2013}.  
Our original motivation was to find an elementary proof that $\rank(\S_m\wr\S_n)=2$.  The argument we found is also applicable to wreath products involving alternating and/or symmetric groups, and indeed to larger classes of wreath products, as described in Lemmas \ref{lem:2gen}--\ref{lem:A}.  

As a first step, we establish a number of two-element generating sets for finite symmetric and alternating groups in Section \ref{sect:SA}, taking as our starting point only the most basic generating sets that have been known for well over a century \cite{Moore1897}.  It might seem somewhat strange that we would need to prove such results, seeing as ``almost all'' pairs of permutations from $\S_n$ generate all of $\S_n$ or $\A_n$, as shown in \cite{Dixon1969}.  See also \cite{Piccard1939} and \cite{IZ1995}; the latter shows that for $n\not=4$, any non-trivial permutation from $\S_n$ is part of a two-element generating set for~$\S_n$.\footnote{This property is known as $\frac32$-generation.  It is known that every finite simple group is $\frac32$-generated \cite{GK2000}.  Although our main result shows that wreath products of finite symmetric or alternating groups are 2-generated, they are not $\frac32$-generated, apart from trivially small exceptions.  For example, if $S$ is not cyclic, then any non-trivial pair $(\ba;\id_n)$ cannot be part of a two-element generating set for $G\wr S$.}  However, we need generating sets with elements of specific orders and with specific fixed points; see especially Lemma~\ref{lem:main}.  In Section \ref{sect:W}, we combine the main results of Section \ref{sect:SA} to prove a number of general (albeit technical) results concerning generation of wreath products, and use these to complete the proof of Theorem~\ref{thm:main}.  In Section \ref{sect:more} we discuss the challenges of extending Theorem \ref{thm:main} to wreath products with more than two factors.  We thank the referee for their helpful comments.

\section{Generating symmetric and alternating groups}\label{sect:SA}

If $a$ is an element of a group, $\ord(a)$ denotes the order of $a$: i.e., the least positive integer $n$ such that $a^n=1$ (the identity of the group) if such an $n$ exists, or else $\ord(a)=\infty$.

\begin{lemma}\label{lem:coprime}
Let $a$ and $b$ be elements of (possibly distinct) groups, and suppose $\ord(a)$ and $\ord(b)$ are finite and coprime.  Suppose also that $r\geq1$ is coprime to $\ord(a)$.  Then there exists $k\geq1$ such that $a^{kr}=a$ and~$b^k=1$.
\end{lemma}

\pf
Write $p=\ord(a)$ and $q=\ord(b)$.  By the Chinese Remainder Theorem, $kr\equiv 1\pmod p$ and $k\equiv0\pmod q$ for some $k\geq1$.
\epf

We now prove some results concerning generators for $\S_n$ and $\A_n$; some of the earlier ones are known (see for example \cite{Conrad,BH2020,Moore1897}), but simple proofs are included for convenience.  In the proofs, for $m\leq n$ we identify~$\S_m$ with the subgroup of $\S_n$ consisting of the permutations fixing $\{m+1,\ldots,n\}$ pointwise.

\begin{lemma}\label{lem:23}
\ben
\item \label{23_1} If $n\geq2$, then $\S_n=\big\la(1,2),(2,3),\ldots,(n-1,n)\big\ra$.
\item \label{23_2} If $n\geq3$, then $\A_n=\big\la(1,2,3),(2,3,4),\ldots,(n-2,n-1,n)\big\ra$
\item[] {\white If $n\geq3$, then} $\phantom{\A_n}= \big\la(1,2,3),(1,2,4),\ldots,(1,2,n)\big\ra$
\item[] {\white If $n\geq3$, then} $\phantom{\A_n}= \big\la(1,2,3),(1,3,4),\ldots,(1,n-1,n)\big\ra$.
\een
\end{lemma}

\pf
Part \ref{23_1} and the first generating set in part \ref{23_2} may be found in \cite{Moore1897}.  (The generators from \ref{23_1} generate $\S_n$ as a Coxeter group \cite{Humphreys1990}.)  Let
\[
U_n=\big\la(1,2,3),(1,2,4),\ldots,(1,2,n)\big\ra \AND V_n=\big\la(1,2,3),(1,3,4),\ldots,(1,n-1,n)\big\ra.
\]
We first show by induction on $n$ that $U_n=\A_n$.  This is clearly true for $n=3$, so suppose $n\geq4$, and assume inductively that $U_{n-1}=\A_{n-1}$.  Then $U_{n-1}$ contains $(1,2,3),(2,3,4),\ldots,(n-3,n-2,n-1)$, and hence so too does $U_n$.  It therefore suffices to show that $(n-2,n-1,n)\in U_n$.  But $(n-2,n-1,n)=(2,1,n-1)(1,2,n)(1,n-2,n-1)$, with $(1,2,n)\in U_n$ and $(2,1,n-1),(1,n-2,n-1)\in\A_{n-1}=U_{n-1}\sub U_n$.

Finally, we can show that $V_n=\A_n$ by showing that $V_n$ contains the generators for $U_n$.  But this follows from $(1,2,i)=(1,2,3)(1,3,4)\cdots(1,i-1,i)$.
\epf


For $f,g\in\S_n$, we write $f^g=g^{-1}fg$.  

\begin{lemma}\label{lem:2n3n}
\ben
\item \label{2n3n_1} If $n\geq2$, then $\S_n=\big\la(1,2),(1,2,3,\ldots,n)\big\ra = \big\la(1,2),(2,3,\ldots,n)\big\ra$.
\item \label{2n3n_2} If $n\geq3$ is odd, then $\A_n=\big\la(1,2,3),(1,2,3,\ldots,n)\big\ra$.
\item \label{2n3n_3} If $n\geq4$ is even, then $\A_n=\big\la(1,2,3),(2,3,\ldots,n)\big\ra$.
\een
\end{lemma}

\pf
\firstpfitem{\ref{2n3n_1}}  This is also contained in \cite{Moore1897}.

\pfitem{\ref{2n3n_2}}  Let $f=(1,2,3)$ and $g=(1,2,3,\ldots,n)$.  Then for any $1\leq i\leq n-2$, $(i,i+1,i+2)=f^{g^{i-1}}$.  The claim then follows from Lemma \ref{lem:23}\ref{23_2}.

\pfitem{\ref{2n3n_3}}  Let $f=(1,2,3)$ and $g=(2,3,\ldots,n)$.  Then for any $1\leq i\leq n-2$, $(1,i+1,i+2)=f^{g^{i-1}}$.  The claim then follows from Lemma \ref{lem:23}\ref{23_2}.
\epf

\begin{lemma}\label{lem:SA}
If $n\geq3$, then
\ben
\item \label{SA_1} $\big\la (1,2,3),(3,4,5,\ldots,n)\big\ra = \begin{cases}
\A_n &\text{if $n$ is odd}\\
\S_n &\text{if $n$ is even,}
\end{cases}$
\item \label{SA_2} $\big\la (1,2,3),(2,3,4,\ldots,n)\big\ra = \begin{cases}
\S_n &\text{if $n$ is odd}\\
\A_n &\text{if $n$ is even.}
\end{cases}$
\een
\end{lemma}

\pf
\firstpfitem{\ref{SA_1}}  Let $S=\la f,g\ra$, where $f=(1,2,3)$ and $g=(3,4,5,\ldots,n)$.  Then $gf=(1,2,3,\ldots,n)$.  As in the proof of Lemma \ref{lem:2n3n}\ref{2n3n_2}, $S$ contains $(i,i+1,i+2)$ for each $1\leq i\leq n-2$.  Thus, $\A_n\sub S$.  If $n$ is odd, then $S\sub\A_n$ (since both $(1,2,3)$ and $(3,4,5,\ldots,n)$ are odd cycles, and hence even permutations), so that $S=\A_n$ in this case.  If $n$ is even, then $S$ contains $\A_n$ and also the odd permutation $(3,4,5,\ldots,n)$, so it follows that $S=\S_n$.

\pfitem{\ref{SA_2}}  The proof is almost identical to the previous part, using the third generating set from Lemma~\ref{lem:23}\ref{23_2}.
\epf

\begin{lemma}\label{lem:SA2}
\ben
\item \label{SA2_1} If $n\geq4$, then $\big\la (1,2)(3,4),(2,3,4,\ldots,n)\big\ra = \begin{cases}
\S_n &\text{if $n\not=5$ is odd}\\
\A_n &\text{if $n$ is even.}
\end{cases}$
\item \label{SA2_2} If $n\geq5$, then $\big\la (1,2)(3,4),(2,4,5,6,\ldots,n)\big\ra = \begin{cases}
\A_n &\text{if $n$ is odd}\\
\S_n &\text{if $n\not=6$ is even.}
\end{cases}$
\een
\end{lemma}

\pf
\firstpfitem{\ref{SA2_1}}  This may be easily verfied with GAP \cite{GAP4} for $n=4,6$, so we assume $n\geq7$ for the rest of the proof.  
Write $f=(1,2)(3,4)$ and $g=(2,3,4,\ldots,n)$, and put $S=\la f,g\ra$.  Then $f^g=(1,3)(4,5)$ and $f^{g^2}=(1,4)(5,6)$.  One may check with GAP that $\A_6=\la f,f^g,f^{g^2}\ra\sub S$.  It follows that $(1,2,3)\in S$.  Since also $(2,3,4,\ldots,n)\in S$, the result then follows from Lemma \ref{lem:SA}\ref{SA_2}.

\pfitem{\ref{SA2_2}}  Again, this may be checked with GAP for $n=5$, so we assume that $n\geq7$.  Write $f=(1,2)(3,4)$ and $g=(2,4,5,6,\ldots,n)$, and put $S=\la f,g\ra$.  This time, $f^g=(1,4)(3,5)$, $f^{g^2}=(1,5)(3,6)$ and $f^{g^3}=(1,6)(3,7)$, and GAP verifies that $\A_7=\la f,f^g,f^{g^2},f^{g^3}\ra\sub S$.  In particular, $S$ contains $(1,2,3)$ and $(2,3,4)$.  Since also $(3,4,5,\ldots,n)=g(2,3,4)\in S$, the result follows from Lemma \ref{lem:SA}\ref{SA_1}.
\epf

\begin{lemma}\label{lem:S0}
If $n=5$ or if $n\geq7$, then $\S_n=\big\la(1,2,3,4),(3,4,5,\ldots,n)\big\ra$.
\end{lemma}

\pf
Again GAP deals with the $n=5$ case, so we assume $n\geq7$.  Write $f=(1,2,3,4)$ and $g=(3,4,5,\ldots,n)$, and put $S=\la f,g\ra$.  GAP shows that $\S_7=\la f,f^g,f^{g^2},f^{g^3}\ra\sub S$.  In particular, $(1,2)$ and $(1,2,3)$ both belong to $S$.  So too therefore does $(1,2,3,\ldots,n)=g(1,2,3)$, and the result then follows from Lemma \ref{lem:2n3n}\ref{2n3n_1}.
\epf

We note that Lemma \ref{lem:SA2}\ref{SA2_1} does not hold for $n=5$, and Lemmas~\ref{lem:SA2}\ref{SA2_2} and \ref{lem:S0} do not hold for $n=6$.
Here is the main technical result of this section:

\begin{lemma}\label{lem:main}
If $S=\S_n$ for some $n\geq4$, or if $S=\A_n$ for some $n\geq5$, then $S=\la f,g\ra$ for some $f,g\in S$ such that
\[
nf=n \COMMA 1g=1 \COMMA \text{$\ord(f)$ is a power of $2$} \COMMA \text{$\ord(g)$ is odd.}
\]
\end{lemma}


\pf
\bit
\item If $S=\S_n$ for even $n\geq4$, then by Lemma \ref{lem:2n3n}\ref{2n3n_1} we may take $f=(1,2)$ and $g=(2,3,4,\ldots,n)$.
\item If $S=\S_n$ for odd $n\geq5$, then by Lemma \ref{lem:S0} we may take $f=(1,2,3,4)$ and $g=(3,4,5,\ldots,n)$.
\item If $S=\A_n$ for even $n\geq6$, then by Lemma \ref{lem:SA2}\ref{SA2_1} we may take $f=(1,2)(3,4)$ and $g=(2,3,4,\ldots,n)$.
\item If $S=\A_n$ for odd $n\geq5$, then by Lemma \ref{lem:SA2}\ref{SA2_2} we may take $f=(1,2)(3,4)$ and $g=(2,4,5,6,\ldots,n)$.  \qedhere
\eit
\epf

We will also need the following simple consequence:

\begin{lemma}\label{lem:main2}
If $G$ is a finite symmetric or alternating group, then $G=\la a,b\ra$ for some $a,b\in G$ such that $\ord(a)$ is odd, and $\ord(b)$ is a power of $2$.
\end{lemma}

\pf
This follows from Lemma \ref{lem:main} unless $G=\S_n$ for some $n\leq3$ or $G=\A_n$ for some $n\leq4$.  These remaining cases are easily checked.  Note that the identity element has order $1=2^0$.
\epf

\section{Generating wreath products}\label{sect:W}

We now turn our attention to wreath products $G\wr S$, where $G$ is a group and $S$ a subgroup of $\S_n$ for some~$n\geq1$.  Although we are primarily interested in the case that $G$ and $S$ are (finite) symmetric or alternating groups, it will be convenient to consider more general situations, as this will allow us to treat multiple cases simultaneously.  We will denote the identity of $\S_n$ (and hence of $S$) by $\id_n$, and the identity of~$G$ by $1$.  We will also denote the $n$-tuple $(1,\ldots,1)$ by $\bone$.  The identity of $G\wr S$ is $(\bone,\id_n)=(1,\ldots,1;\id_n)$.  In the following lemmas, $G$ could be (countably) infinite, and we recall that any generating sets are assumed to be semigroup generating sets.

\begin{lemma}\label{lem:4gen}
Suppose $S=\la f,g\ra$ is a transitive subgroup of $\S_n$ for some $n\geq1$, and $G$ is a group with $G=\la a,b\ra$.  Then for arbitrary $h_1,h_2\in S$, $G\wr S=\la\al,\be,\ga,\de\ra$, where
\[
\al=(\bone,f) \COMMA \be=(\bone,g) \COMMA \ga=(1,\ldots,1,a,1,\ldots,1;h_1) \COMMA \de=(1,\ldots,1,b,1,\ldots,1;h_2),
\]
where the $a$ and $b$ can be in arbitrary coordinates in $\ga$ and $\de$.
\end{lemma}

\pf
Let $P=\la\al,\be,\ga,\de\ra$.  Certainly $P\sub G\wr S$.  Since $S=\la f,g\ra$, $P$ contains every element of the form~$(\bone,h)$, for $h\in S$.

Suppose $a$ is in the $i$th coordinate of $\ga$.  Since $S$ is transitive, there exists $h\in P$ such that $1h=i$.  But then
\[
\ga' := (a,1,\ldots,1;\id_n) = (\bone,h) \cdot \ga\cdot (\bone,h_1^{-1}h^{-1}) \in P.
\]
A similar calculation shows that $\de':=(b,1,\ldots,1;\id_n)\in P$.  
Since $G=\la a,b\ra$, it quickly follows that every element of the form $(c,1,\ldots,1;\id_n)$, for $c\in G$, is contained in $\la\ga',\de'\ra\sub P$.  

By transitivity of $S$ again, $P$ contains every element of the form $(1,\ldots,1,c,1,\ldots,1;\id_n)$, with $c\in G$ in an arbitrary coordinate.  It then quickly follows that $P$ contains every element of the form $(\ba,\id_n)$ for~$\ba\in G^n$.  But then for any $\ba\in G^n$ and $h\in S$, $(\ba,h)=(\ba,\id_n)\cdot(\bone,h)\in P$.  This shows that $G\wr S\sub P$.
\epf

\begin{lemma}\label{lem:2gen}
Suppose $S=\la f,g\ra$ is a transitive subgroup of $\S_n$ for some $n\geq4$, where
\[
nf=n \COMMA 1g=1 \COMMA \text{$\ord(f)$ is a power of $2$} \COMMA \text{$\ord(g)$ is odd.}
\]
Suppose $G=\la a,b\ra$ is a group, where $\ord(a)$ is odd, and $\ord(b)$ is a power of $2$.  Then 
\[
G\wr S=\la\al,\be\ra, \qquad\text{where}\qquad \al = (1,\ldots,1,a;f) \AND \be = (b,1,\ldots,1;g).
\]
\end{lemma}

\pf
Let $P=\la\al,\be\ra$.  Certainly $P\sub G\wr S$.  Since $\ord(a)$ and $\ord(f)$ are coprime, it follows from Lemma~\ref{lem:coprime} (with $r=1$) that $a^k=a$ and $f^k=\id_n$ for some $k\geq1$.  Similarly, $b^l=b$ and $g^l=\id_n$ for some $l\geq1$.  But then, using $nf=n$ and $1g=1$, we see that
\[
\ga:=(1,\ldots,1,a;\id_n)=(1,\ldots,1,a^k;f^k)=\al^k\in P \ANDSIM \de: = (b,1,\ldots,1;\id_n)=\be^l\in P.
\]
Since $\ord(\ga)=\ord(a)$ and $\ord(\de)=\ord(b)$ are finite, $\ga^{-1}$ and $\de^{-1}$ both belong to $P$.  So too therefore do 
\[
\al' := \ga^{-1}\al = (\bone,f) \AND \be' := \de^{-1}\be = (\bone,g).
\]
It follows from Lemma \ref{lem:4gen} that $G\wr S=\la\al',\be',\ga,\de\ra\sub P$.
\epf

We will use Lemma \ref{lem:2gen} to establish Theorem \ref{thm:main} in the generic case in which $S$ has large enough degree.  For the cases in which $S$ has small degree, we require the following two results.

\begin{lemma}\label{lem:S}
Suppose $G=\la a,b\ra$ is a group, where $\ord(a)$ is odd, and $\ord(b)$ is a power of $2$.  Then
\ben
\item \label{S1}  $G\wr\S_2=\la\al,\be\ra$, where $\al=(ab^{-1},1;\id_2)$ and $\be=(1,b;(1,2))$,
\item \label{S2}  $G\wr\S_3=\la\al,\be\ra$, where $\al=(a,b,1;(2,3))$ and $\be=(1,1,1;(1,2))$.
\een
\end{lemma}

\pf
For both parts we write $P=\la\al,\be\ra$.  As above, we just need to show that $G\wr\S_n\sub P$ ($n=2,3$).

\pfitem{\ref{S1}}  First note that $\al\be^2=(a,b;\id_2)$.  Let $k,l\geq1$ be such that $a^k=a$, $b^l=b$ and $a^l=b^k=1$ (cf.~Lemma~\ref{lem:coprime}).  Then~$P$ contains
\[
\ga:=(\al\be^2)^k=(a^k,b^k;\id_2)=(a,1;\id_2) \AND \de := (\al\be^2)^l=(a^l,b^l;\id_2)=(1,b;\id_2).
\]
Since $\ord(\de)=\ord(b)$ is finite, $P$ contains $\de^{-1}$, and hence also $\al':=\de^{-1}\be=(1,1;(1,2))$.  Lemma~\ref{lem:4gen} then gives $G\wr\S_2=\la\al',\ga,\de\ra\sub P$.

\pfitem{\ref{S2}}  Here we note that ${\al^2=(a^2,b,b;\id_3)}$.  By Lemma \ref{lem:coprime} (with $r=2$), $a^{2k}=a$ and $b^k=1$ for some ${k\geq1}$.  But then~$P$ contains $\ga:=\al^{2k}=(a^{2k},b^k,b^k;\id_3)=(a,1,1;\id_3)$, and hence also ${\de:=\ga^{-1}\al=(1,b,1;(2,3))}$.  Then with $\al':=\de\be\de^{-1}=(1,1,1;(1,3))$, Lemma \ref{lem:4gen} gives $G\wr\S_3=\la\al',\be,\ga,\de\ra\sub P$.
\epf

\begin{lemma}\label{lem:A}
Suppose $G=\la a,b\ra$ is a group, where $\ord(a)$ is odd, and $\ord(b)$ is a power of $2$.  Then
\ben
\item \label{A1}  $G\wr\A_2=\la\al,\be\ra$, where $\al=(a,b;\id_2)$ and $\be=(b,a;\id_2)$,
\item \label{A2}  $G\wr\A_3=\la\al,\be\ra$, where $\al=(ab^{-1},1,1;\id_3)$ and $\be=(b,1,1;(1,2,3))$,
\item \label{A3}  $G\wr\A_4=\la\al,\be\ra$, where $\al=(a,1,1,b;(1,2,3))$ and $\be=(1,1,1,1;(2,3,4))$.
\een
\end{lemma}

\pf
For all parts we write $P=\la\al,\be\ra$.  As usual, we must show that $G\wr\A_n\sub P$ ($n=2,3,4$).

\pfitem{\ref{A1}}  Let $k,l\geq1$ be such that $a^k=a$, $b^l=b$ and $a^l=b^k=1$ (cf.~Lemma \ref{lem:coprime}).  Then the claim follows from
\[
\al^k=(a,1;\id_2) \COMMA
\be^l=(b,1;\id_2) \COMMA
\be^k=(1,a;\id_2) \COMMA
\al^l=(1,b;\id_2) .
\]

\pfitem{\ref{A2}}  First note that $\al\be^3=(a,b,b;\id_3)$.  Let $k\geq1$ be such that $a^k=a$ and $b^k=1$ (cf.~Lemma \ref{lem:coprime}).  Then~$P$ contains $\ga:=(\al\be^3)^k=(a^k,b^k,b^k;\id_3)=(a,1,1;\id_3)$, and hence also $\al':=\ga^{-1}\al\be=(1,1,1;(1,2,3))$.  Lemma \ref{lem:4gen} then gives $G\wr\A_3=\la\al',\be,\ga\ra\sub P$.

\pfitem{\ref{A3}}  Now, $\al^3=(a,a,a,b^3;\id_4)$.  By Lemma \ref{lem:coprime} (with $r=3$, and with the roles of $a$ and $b$ swapped), $a^k=1$ and ${b^{3k}=b}$ for some $k\geq1$.  So~$P$ contains $\ga:=\al^{3k}=(a^k,a^k,a^k,b^{3k};\id_4)=(1,1,1,b;\id_4)$, and hence also $\de:=\ga^{-1}\al=(a,1,1,1;(1,2,3))$, and $\al':=\de^{-1}\be\de=(1,1,1,1;(1,4,3))$.  Since $\big\la(1,4,3),(2,3,4)\big\ra=\A_4$, Lemma~\ref{lem:4gen} then gives $G\wr\A_4=\la\al',\be,\ga,\de\ra\sub P$.
\epf

We may now tie together the loose ends.

\pf[\bf Proof of Theorem \ref{thm:main}]
Suppose $G$ and $S$ are arbitrary finite symmetric or alternating groups.  

We begin with the second statement.  First, if either of \ref{cyc1} or \ref{cyc2} holds, then $G\wr S$ is clearly cyclic.  For the converse, we make two observations:
\bit
\item If either $G$ or $S$ is not cyclic, then it is non-abelian (recall that $G$ and $S$ are symmetric or alternating groups), in which case $G\wr S$ is also non-abelian, and therefore not cyclic.  
\item If $G$ and $S$ are both non-trivial, say with non-identity elements $f\in G$ and $g\in S$, then the elements $(1,\ldots,1,f,1,\ldots,1;\id_n)$ and $(1,\ldots,1;g)$ do not commute (here the $f$ is in position $i$ for some $i\not=ig$), so again $G\wr S$ is not cyclic.
\eit
The above two points show that if $G\wr S$ is cyclic, then $G$ and $S$ are both cyclic, and at most one of them is non-trivial.  This still includes the cases in which $G$ is non-trivial and cyclic (i.e., one of $\S_2$ or $\A_3$), and $S=\A_2$, which are not listed in the theorem.  But in these cases, $G\wr S\cong G^2$ is non-cyclic.

We now prove the first statement: i.e., that $G\wr S$ is 2-generated.  By Lemma \ref{lem:main2}, $G=\la a,b\ra$ for some $a,b$ with $\ord(a)$ odd, and $\ord(b)$ a power of 2.  The cases in which $S=\S_n$ ($n=2,3$) or $S=\A_n$ ($n=2,3,4$) are covered by Lemmas \ref{lem:S} and \ref{lem:A}.  The case in which $S=\S_1=\A_1$ is clear.  
If $S=\S_n$ ($n\geq4$) or $S=\A_n$ ($n\geq5$), then Lemma \ref{lem:main} tells us that Lemma~\ref{lem:2gen} applies to $S$, and the conclusion of that lemma completes the proof.
\epf

\section{Wreath products with more than two factors}\label{sect:more}

We believe it would be interesting to investigate minimal generation of iterated wreath products of the form
\[
G_1\wr G_2\wr G_3\wr\cdots\wr G_k := (\cdots((G_1\wr G_2)\wr G_3)\wr\cdots)\wr G_k,
\]
where each factor is a finite symmetric or alternating group.  For example, Table \ref{tab:d} gives values of $d(G_1\wr G_2\wr G_3)$ where the $G_i$ are symmetric or alternating groups of relatively small degree, computed using GAP \cite{GAP4}.  Note that the rows are indexed by $G_1\wr G_2$, and the columns by $G_3$; thus, for example, $d(\S_3\wr\S_3\wr\A_2)=4$ and $d(\S_3\wr\A_3\wr\A_2)=2$.
%
Even this limited data shows that a formula for $d(G_1\wr G_2\wr G_3\wr\cdots\wr G_k)$ is likely to be rather complicated for $k\geq3$.  On the other hand, the proof of Lemma \ref{lem:4gen} may easily be modified to show that $d(G\wr S) \leq d(G) + d(S)$ for any group $G$ and any transitive subgroup $S$ of $\S_n$.  
Thus, if $G_1,\ldots,G_k$ are non-trivial finite symmetric or alternating groups (so we are excluding the non-transitive $\A_2$), then
\[
d(G_1\wr G_2\wr G_3\wr\cdots\wr G_k) \leq d(G_1\wr G_2) + d(G_3) + \cdots + d(G_k) \leq 2(k-1).
\]
We note that this inequality does not hold in general if we allow $\A_2$ factors.  For example, since $G\wr\A_2\cong G^2$ for any group $G$, it follows that $\S_2\wr\A_2\wr\A_2\wr\A_2\cong\S_2^8$, which is minimally 8-generated.

\begin{table}[ht]
\begin{center}
\begin{tabular}{l|cccccc}
 & $\A_2$ & $\A_3$ & $\A_4$ & $\S_2$ & $\S_3$ & $\S_4$ \\
\hline
$\A_2\wr\A_2$ & 1 & 1 & 2 & 1 & 2 & 2 \\
$\A_2\wr\A_3$ & 2 & 2 & 2 & 2 & 2 & 2 \\ 
$\A_2\wr\S_2$ & 2 & 2 & 2 & 2 & 2 & 2 \\ 
$\A_2\wr\S_3$ & 2 & 2 & 2 & 2 & 2 & 2 \\ 
$\A_3\wr\A_2$ & 4 & 3 & 3 & 3 & 2 & 2 \\ 
$\A_3\wr\A_3$ & 4 & 3 & 3 & 3 & 2 & 2 \\ 
$\A_3\wr\S_2$ & 2 & 2 & 2 & 2 & 2 & 2 \\ 
$\A_3\wr\S_3$ & 2 & 2 & 2 & 2 & 2 & 2 \\ 
$\S_2\wr\A_2$ & 4 & 3 & 2 & 3 & 3 & 3 \\ 
$\S_2\wr\A_3$ & 2 & 2 & 2 & 2 & 2 & 2 \\ 
$\S_2\wr\S_2$ & 4 & 3 & 2 & 3 & 3 & 3 \\ 
$\S_2\wr\S_3$ & 4 & 3 & 2 & 3 & 3 & 3 \\ 
$\S_3\wr\A_2$ & 4 & 3 & 2 & 3 & 3 & 3 \\ 
$\S_3\wr\A_3$ & 2 & 2 & 2 & 2 & 2 & 2 \\ 
$\S_3\wr\S_2$ & 4 & 3 & 2 & 3 & 3 & 3 \\ 
$\S_3\wr\S_3$ & 4 & 3 & 2 & 3 & 3 & 3 \\ 
\end{tabular}
\end{center}
\vspace{-.5cm}
\caption{Values of $d(G_1\wr G_2\wr G_3)$, where the $G_i$ are symmetric or alternating groups of small degree.}
\label{tab:d}
\end{table}

\footnotesize
\def\bibspacing{-1.1pt}
\bibliography{biblio}
\bibliographystyle{abbrv}
\end{document}